\documentclass[12pt]{article}
\usepackage{latexsym}
\usepackage{amsfonts}
\def\Bbb{\mathbb}

\def\eea{\end{eqnarray*}}

\newtheorem{main}{Theorem}

\newtheorem{thm}{Theorem}[section]
\newtheorem{prop}[thm]{Proposition}
\newtheorem{cor}[thm]{Corollary}
\newtheorem{lem}[thm]{Lemma}
\newtheorem{defn}[thm]{Definition}
\newtheorem{conjecture}[thm]{Conjecture}

\newenvironment{proof}{\medskip \noindent
{\bf Proof.}}{\hfill \rule{.5em}{1em}
\\}

\begin{document}
\sloppy
\title{Hyperbolic Manifolds, Harmonic Forms, \\
and  Seiberg-Witten Invariants}

\author{Claude LeBrun\thanks{Supported 
in part by  NSF grant DMS-0072591.} 
\\ 
SUNY Stony
 Brook 
  }

\date{}

\maketitle

\begin{abstract}
New estimates are derived concerning the behavior
of self-dual hamonic $2$-forms on a 
compact Riemannian $4$-manifold with
non-trivial Seiberg-Witten invariants. 
Applications include a vanishing theorem for
certain Seiberg-Witten invariants on 
 compact $4$-manifolds of constant 
negative sectional  curvature. 
 \end{abstract}

\section{Introduction}

Seiberg-Witten theory gives rise to  a powerful interplay 
between the geometry and topology of smooth $4$-manifolds. 
For example, a remarkable theorem of Taubes \cite{taubes} asserts that 
any symplectic  $4$-manifold with $b^{+}\geq 2$ has
a non-zero Seiberg-Witten invariant.
On the other hand, if a  $4$-manifold with $b^{+}\geq 2$ 
admits a metric of positive scalar curvature, Witten 
\cite{witten} observed
that its 
 Seiberg-Witten invariants must all vanish. 
The existence of a metric satisfying a suitable
curvature condition may thus be sufficient to 
rule out the existence of a symplectic structure
on a given smooth,  compact $4$-manifold.

Despite this, there are many   ``naturally occurring''  classes of
 $4$-manifolds for which we do not yet know whether 
 the Seiberg-Witten invariants all vanish. In particular, 
 the following still appears to be open: 

\begin{conjecture}
Let $M^{4}= {\mathcal H}^{4}/\Gamma$ be a compact hyperbolic
$4$-manifold. Then all the Seiberg-Witten invariants 
of $M$ vanish. In particular, $M$ does not admit symplectic 
structures. 
\end{conjecture}

By contrast, the complex-hyperbolic $4$-manifolds, which
 by definition are compact quotients of the
unit ball in ${\mathbb C}^{2}$,  are all K\"ahler manifolds,
and so carry  symplectic structures compatible with their
standard  orientations. One might expect, however,
for the situation to be quite different 
regarding the {\em non-standard} orientation of 
 such a manifold:

\begin{conjecture}
Let $M^{4}= \overline{{\mathbb C}{\mathcal H}_{2}/\Gamma}$ be a
{\em reverse-oriented} compact complex-hyperbolic
$4$-manifold. Then all the Seiberg-Witten invariants 
of $M$ vanish for the fixed orientation. In particular, 
there is no symplectic structure on $M$ compatible
with the non-complex orientation. 
\end{conjecture}

In fact, both of these speculations may be interpreted as 
special cases of a more general conjecture. 
Recall that the $2$-forms on an oriented $4$-manifold
decompose as 
$$\Lambda^{2}= \Lambda^{+}\oplus \Lambda^{-},$$ 
where $\Lambda^{\pm}$ is the $(\pm 1)$ eigenspace
of  Hodge  star operator $\star$. Thinking of the 
curvature tensor $\mathcal R$ as a linear map
$\Lambda^{2}\to \Lambda^{2}$, we thus \cite{st} get 
a decomposition 
$$
{\mathcal R}=
\left(
\mbox{
\begin{tabular}{c|c}
&\\
$W_++\frac{s}{12}$&$\stackrel{\circ}{r}$\\ &\\
\cline{1-2}&\\
$\stackrel{\circ}{r}$ & $W_{-}+\frac{s}{12}$\\&\\
\end{tabular}
} \right) 
$$
into irreducible pieces. 
Here the
self-dual and anti-self-dual Weyl curvatures 
$W_\pm$ are the trace-free pieces of the appropriate blocks.
The scalar curvature  $s$ is understood to act by scalar multiplication,
whereas the trace-free part 
$\stackrel{\circ}{r}=r-\frac{s}{4}g$  of the Ricci curvature
acts on 2-forms by
$$\psi_{ab} \mapsto ~ 
 \stackrel{\circ}{r}_{ac}{\psi^c}_{b}-
\stackrel{\circ}{r}_{bc}{\psi^c}_{a}.$$
An oriented $4$-manifold is said to be 
{\em self-dual} if $W_{-}\equiv 0$, or {\em anti-self-dual}
if $W_{+}\equiv 0$. 
 An orientable real-hyperbolic $4$-manifold is 
{\em both} self-dual and anti-self-dual, which
is just another way of saying that any such  manifold
is locally conformally flat. On the other hand, a complex-hyperbolic 
$4$-manifold is self-dual with respect to the orientation
determined by the complex structure; 
thus, since reversing the orientation of 
a $4$-manifold interchanges $W_+$ and $W_-$,
a complex-hyperbolic 
$4$-manifold is 
 anti-self-dual with respect to its {\em non-complex}
orientation. Also notice that both real- and complex-hyperbolic
$4$-manifolds are {\em Einstein} --- i.e. they satisfy 
$\stackrel{\circ}{r}\equiv 0$. Thus the above conjectures
might be viewed as simply special cases of the following:

\begin{conjecture}
Let $(M^{4},g)$ be a compact anti-self-dual Einstein manifold
with negative scalar curvature. Then all the Seiberg-Witten invariants 
 vanish for the fixed orientation of $M$.
\end{conjecture}

This paper will present some tantalizing, albeit inconclusive, 
 evidence in
favor of these conjectures. 
To this end, let us first draw the reader's attention to a 
beautiful recent result  of Armstrong \cite{arm1}  asserting 
 that anti-self-dual Einstein spaces of negative scalar curvature 
 never admit
non-trivial self-dual harmonic $2$-forms of 
constant length. On the other hand, 
this article will prove that if such a space 
has a non-trivial
Seiberg-Witten invariant, it necessarily admits a self-dual  harmonic 2-form
whose length is ``nearly constant,'' by  two different
quantitative measures. 
A quantitative sharpening of Armstrong's result might
therefore provide exactly the tool needed to  prove some version
of the above conjectures.

\section{Harmonic $2$-Forms}
\label{harm} 

 In this section, we introduce 
two different invariants which  offer quantitative obstructions
to the existence of  non-trivial
 self-dual
 harmonic $2$-forms  of constant length
 on a given $4$-manifold.

Our first invariant is simply the minimal angle 
between the point-wise norm of the form and the constant
$1$, considered as vectors in the Hilbert space $L^{2}$:

\begin{defn}\label{angl}
Let $(M,g)$ be a compact, oriented  Riemannian $4$-manifold
with $b^{+}(M) \geq 1$.
Let 
 $${\cal H}^{+}_{g} = \{ \phi \in {\mathcal E}^{2}(M)~|~ \phi = \star 
 \phi , ~~ d\phi =0 \}$$
 be the space of self-dual harmonic 2-forms on $(M,g)$, so that 
  $\dim {\cal H}^{+}_{g}  = b^{+}(M)  > 0$. 
 We define
$$\theta (M,g) = \min_{\phi \in ({\mathcal H}^{+}_{g}-0)}   
\cos^{-1}
\left( 
\frac{\int_{M}|\phi |d\mu_{g}}{V^{1/2}
\left( \int_{M}|\phi |^{2}d\mu\right)^{1/2}}
\right)$$
where  $V=\int_{M}1~d\mu_{g}$ is the total volume
of $(M,g)$.
\end{defn}

Our second invariant is rather more subtle, 
and is best understood
in the context of the following observation:

\begin{prop}
Let $(M,g)$ be any  compact oriented Riemannian $4$-manifold, and 
let $\phi$ be any self-dual harmonic $2$-form. Then the 
 function  $f = \sqrt{|\phi |}$
 belongs to the Sobolev space $L^{2}_{1}$, and satisfies 
 \begin{equation}
	\int_{M} |df|^{2}d\mu \leq \int_{M} 
\left( \sqrt{\frac{2}{3}}|W_{+}|-\frac{s}{6} \right) |\phi | ~ d\mu  ,
	\label{super}
\end{equation}
\end{prop}
\begin{proof}
Since $\phi$ is smooth by elliptic regularity, $f$ is certainly 
continuous, and thus belongs to $L^{2}$. We therefore just need to 
show is that $|df|$
 belongs to $L^{2}$, and satisfies (\ref{super}). 
 
 We may assume henceforth that $\phi\not\equiv 0$, since otherwise there
 is nothing to prove. It then follows that the {\em nodal set}
 where $\phi$ vanishes is of measure zero;   indeed \cite{baer}, its  
 Hausdorff dimension is $\leq 2$. The function $|df|^{2}$ is 
 smooth outside this nodal set, and our objective is just 
 to show that its integral over the complement of the 
 nodal set is finite. 

To this end, recall  that the harmonicity of $\phi$ implies
the Weitzenb\"ock formula \cite{bourg} 
\begin{eqnarray*}
	0 & = & \frac{1}{2}\Delta |\phi |^{2} + |\nabla \phi |^{2} - 2 
W_{+}(\phi , \phi ) + \frac{s}{3}|\phi |^{2}  \\
	 & \geq  &  \frac{1}{2}\Delta |\phi |^{2} + |\nabla \phi |^{2}
 +\left( \frac{s}{3}-2\sqrt{\frac{2}{3}}|W_{+}|\right)
|\phi |^{2}
\end{eqnarray*}
On the open set defined by  $\phi\neq 0$, we therefore have
\begin{eqnarray*}
	0 & \geq  & \frac{\Delta |\phi |^{2}}{2 |\phi |}+
	\frac{|\nabla \phi |^{2}}{|\phi |}
+\left(\frac{s}{3}- 2\sqrt{\frac{2}{3}}|W_{+}| \right) |\phi |  \\
	 & = & \Delta |\phi |-
	 \frac{|~d |\phi | ~ |^{2}}{ |\phi |}+\frac{|\nabla \phi |^{2}}{|\phi|}
+\left(\frac{s}{3}- 2\sqrt{\frac{2}{3}}|W_{+}| \right) |\phi |
\end{eqnarray*}
But since $\phi$ is harmonic, we therefore have the refined
Kato inequality \cite{herz}
$$
|\nabla  \phi |^{2}\geq \frac{3}{2} |~d |\phi | ~|^{2} ,
$$
so the above yields 
$$0 \geq \Delta |\phi | + \frac{1}{2}  \frac{|~d |\phi | ~ |^{2}}{ |\phi |}
+\left(\frac{s}{3}- 2\sqrt{\frac{2}{3}}|W_{+}| \right) |\phi | , 
$$
and we thus have 
$$( 2\sqrt{\frac{2}{3}}|W_{+}|-\frac{s}{3} ) |\phi | \geq \Delta 
|\phi | + 
 2|df|^{2}
$$
wherever $\phi \neq 0$. 

Now let $F=|\phi | = f^{2}$, and let $\varepsilon^{2}$ be any 
positive regular value of the smooth function $F^{2}$. Let
$M_{\varepsilon}$ be the set where 
$F \leq \varepsilon$, and observe that, since $\Delta$
is the {\em positive} Laplacian, Stokes' theorem tells
us that 
$$
\int_{M-M_{\varepsilon}}\Delta |\phi | ~d\mu	  =  -  
\int_{M-M_{\varepsilon}}d\star dF
	  =  \int_{\partial M_{\varepsilon}}\star dF   
	  >  0.$$
Thus 
\begin{eqnarray*}
	\int_{M-M_{\varepsilon}} 
( 2\sqrt{\frac{2}{3}}|W_{+}|-\frac{s}{3} ) |\phi | ~ d\mu  & \geq  &
 \int_{M-M_{\varepsilon}}\Delta |\phi | ~d\mu
+ 2\int_{M-M_{\varepsilon}}|df|^{2}d\mu\\
	 & \geq &  2\int_{M-M_{\varepsilon}}|df|^{2}d\mu 
\end{eqnarray*}
On the other hand, $\bigcap_{\varepsilon > 0}M_{\varepsilon}$ 
is the nodal set, which 
has measure zero in $(M,d\mu)$.
Thus
\begin{eqnarray*}
	\int_{M} |df|^{2}d\mu & = & \limsup_{\varepsilon  \searrow 0}
\int_{M-M_{\varepsilon}} |df|^{2}d\mu  \\
	 &\leq  & \limsup_{\varepsilon  \searrow 0} 
	 \int_{M-M_{\varepsilon}} 
( \sqrt{\frac{2}{3}}|W_{+}|-\frac{s}{6} ) |\phi | ~ d\mu \\
	 & = & \int_{M} 
( \sqrt{\frac{2}{3}}|W_{+}|-\frac{s}{6} ) |\phi | ~ d\mu < \infty ,
\end{eqnarray*}
so that $|df|$ is an $L^{2}$ function, and satisfies the promised 
estimate. 
\end{proof}

In particular,  the following invariant  is {\em a priori} finite:

\begin{defn}
\label{quot}
Let $(M,g)$ be a compact oriented  Riemannian $4$-manifold
with $b^{+}(M) \geq 1$. 
 We then define
 $$\nu  (M,g) := \min_{{\phi}\in ({\cal H}^{+}_{g}-0)}
 \frac{\int \left |d\sqrt{|\phi |}\right |^{2} d\mu_{g}}{\int  
 |\phi | d\mu_{g}}$$
\end{defn}

  One might do well to compare this definition  with that of 
  the spectral	invariant 
  $$\lambda_{1}(M,g) =	\inf \left\{  
  \frac{\int  |df |^{2} d\mu_{g}}{\int	 
 	|f|^{2}	d\mu_{g}} ~\left| ~	f\in L^{2}_{1}(M,g)	,
 	 f\not\equiv 0,	\int f d\mu	= 0	\right.\right\} .
 	$$
 	By contrast, 
 	$$\nu  (M,g) = \inf	\left\{	 
  \frac{\int  |df |^{2} d\mu_{g}}{\int	 
 	|f|^{2}	d\mu_{g}} ~\left| ~		\exists	\phi \in 
 	({\cal H}^{+}_{g}-0) \mbox{  s.t. }  f =	\sqrt{|\phi	|}\right.\right\} .
 	$$
 	Despite the analogy, however, these invariants
 	would seem to have little to say about each other. 
 	On one hand, $\nu$ is defined in terms of a finite-dimensional
 	family of functions $f$; on the other hand, 
 	 these functions  are
 	not orthogonal to the constants!
 %

The invariant $\nu$ is not scale-invariant; if the
metric $g$ is replaced by $cg$, where $c$ is a
positive constant, $\nu$ gets replaced by $c^{-1}\nu$;
thus $\nu$  rescales in anology to the scalar curvature $s$ 
or  the first eigenvalue 
$\lambda_{1}$ of the Laplacian. 
 As we are primarily concerned here with manifolds of constant negative scalar 
curvature, we shall  simplify the 
statements of many of our results by assuming the metric in question
satisfies  
$s\equiv -12$,  as does the standard  $K=-1$ metric 
 on a real-hyperbolic $4$-manifold. For example:

\begin{prop}
Let $(M,g)$ be an compact  anti-self-dual Einstein manifold with $s=-12$.
Then $\theta (M,g) \neq 0$, and $\nu (M,g) \in (0, 2]$. 
\end{prop}
\begin{proof} Since $s\equiv -12$ and $W_{+}\equiv 0$ by assumption, 
 the inequality
(\ref{super}) tells us that 
$$\int |d\sqrt{|\phi |}|^{2}d\mu \leq 2\int |\phi |d\mu$$
for any  $\phi \in {\cal H}^{+}-0$. This shows that 
$\nu (M,g ) \leq 2$, as claimed.

If, on the other hand,  either 
$\theta (M,g)$ or $\nu (M,g )$ were zero, $(M,g)$ would admit  
a self-dual harmonic 2-form $\phi$ of constant, non-zero length. 
Multiplying this form by  a constant
would then give us a harmonic self-dual $2$-form $\omega$ of 
point-wise norm  $\equiv \sqrt{2}$. This symplectic form
would then be expressible as $\omega (\cdot , \cdot ) = 
g(J\cdot , \cdot )$ for a unique almost-complex structure
$J$ on $M$, making $(M,g,J,\omega )$ an {\em almost-K\"ahler
manifold}. However, Armstrong \cite{arm1} has shown that 
compact anti-self-dual almost-K\"ahler Einstein manifolds
with $s < 0$ 
do not exist. By contradiction, 
$\nu (M, g)$ and $\theta (M,g)$  must therefore be positive. 
\end{proof}

By contrast,   $\nu=0$  for 
any $4$-dimensional K\"ahler manifold, since the K\"ahler form is
a self-dual harmonic $2$-form of 
constant length. Nonetheless, K\"ahler
manifolds give us an instructive set of examples
when we examine the harmonic  
forms which are {\em orthogonal} to the K\"ahler form:

\begin{prop}\label{pomp}
Let $(M^{4},g)$ be any compact K\"ahler manifold with 
$s=-12$, and let $\phi\not\equiv 0$ be any self-dual harmonic $2$-form
which is $L^{2}$-orthogonal to the K\"ahler form 
$\omega$. Then 
$$
\frac{\int_{M}|d\sqrt{|\phi |}|^{2}d\mu}{\int_{M} |\phi |d\mu}=
\frac{3}{2}.
$$
\end{prop}
\begin{proof}
On a K\"ahler manifold of real dimension $4$, 
$$\Lambda^{+}\otimes {\mathbb C}={\mathbb C}\omega \oplus 
\Lambda^{2,0}\oplus \Lambda^{0,2},$$
and any real self-dual harmonic $2$-form can 
be uniquely written  as 
$$\phi = a \omega + \varphi + \bar{\varphi}, $$
where $a$ is a real constant and $\varphi$ is a holomorphic
$(2,0)$-form. If $\phi$ is $L^{2}$-orthogonal to $\omega$,
$a=0$, and $|\phi |^{2}= 2|\varphi |^{2}$. 
But the Ricci form $\rho$ of $(M,\omega )$ is given by
$$\rho = 
i\partial \bar{\partial} \log|\varphi |^{2}= 
i\partial \bar{\partial} \log |\phi |^{2}$$
away from the nodal set of $\phi$, 
and taking the trace against $\omega$ yields 
$$s = -\Delta \log |\phi |^{2} = -4\Delta \log f,
$$ 
where $f=\sqrt{|\phi |}$.  
Thus 
$$ -\frac{s}{4}f^{2}=f\Delta f + |df|^{2}
$$
away from the nodal set. Again setting  
$F=|\phi|=f^{2}$, and letting $M_{\varepsilon}$ be the set where
$F 
<\varepsilon \}$, for $\varepsilon >0$ any  regular value 
of $F$, we have 
$$-\int_{M-M_{\varepsilon}}\frac{s}{4}f^{2}d\mu =
2\int_{M-M_{\varepsilon}}|df|^{2}d\mu + \frac{1}{2}
\int_{\partial M_{\varepsilon}} \star dF.$$
On the other hand
$$\left|\int_{\partial M_{\varepsilon}} \star dF\right| < C\varepsilon$$
because $|dF|= |~d|\phi||\leq |\nabla \phi|$ by the Kato inequality, 
and the $3$-dimensional volume of the hypersurface $\partial M_{\varepsilon}$
is   $< \mbox{const} \cdot \varepsilon$ by
 the Weierstrass preparation theorem for
holomorphic functions.   Taking the limit as 
$\varepsilon \searrow 0$, we thus obtain
$$-\int_{M} \frac{s}{4} f^{2}d\mu = 2\int_{M}|df|^{2}d\mu.$$
Setting $s=-12$, this tells us that 
$$\frac{\int_{M}|df|^{2}d\mu}{\int_{M} f^{2}d\mu}= -\frac{s}{8}= \frac{3}{2},$$
as claimed.
\end{proof}

Thus, it does not seem unreasonable to hope that   the value of $\nu$
may turn out to be of the order of $1$ for  many manifolds
with $s=-12$. 
With this  in mind, 
we now  state our main result:

\begin{main} \label{wild}
Let $(M,g)$ be a compact  anti-self-dual Einstein manifold with $s=-12$
and $b^{+}(M) \geq 2$.
If $\nu (M, g) \geq 2-\sqrt{3}\approx 0.268$,
all the  
Seiberg-Witten invariants of  the oriented $4$-manifold
$M$ must vanish. 
\end{main}

\section{Seiberg-Witten Estimates}
\label{key}

In this section, we will derive a new set of 
Seiberg-Witten estimates by combining
ideas previously used in \cite{lcp2} and \cite{lric}.
While these estimates suffice to imply
the vanishing results contained in the last section
of the paper, they also have an intrinsic interest of their own, as
well as other ramifications which 
would seem to be worthy of  exploration.

Let $M$ be a smooth compact oriented connected $4$-manifold with 
$b^{+}(M)\geq 2$, and let 
$\mathfrak c$ be any spin$^{c}$ structure on $M$. 
For any Riemannian metric $g$ on $M$, 
we then have rank-$2$ complex Hermitian vector bundles
${\mathbb V}_{\pm}\to M$ 
which formally satisfy 
$${\mathbb V}_{\pm}={\Bbb S}_{\pm}\otimes L^{1/2},$$
where ${\Bbb S}_{\pm}$ are the locally-defined left- and right-handed 
spinor bundles of $g$, and $L=\wedge^{2}{\mathbb V}_{\pm}$ is a globally defined
Hermitian 
line bundle. As a convenient abuse, we will  use  
$c_{1}$ to denote the image of $c_{1}(L)=c_{1}({\mathbb V}_{\pm})
\in H^{2}(M,{\mathbb Z})$ 
in the real cohomology $H^{2}(M,{\mathbb R})$, and refer $c_{1}$ as the 
first   Chern class of $\mathfrak c$. Now the Hodge theorem
tells us that $H^{2}(M,{\mathbb R})$ can be identified
with the space of harmonic $2$-forms ${\mathcal H}^{2}_{g}$
on $(M,g)$; and the latter splits as the direct sum
$${\mathcal H}^{2}_{g}= {\mathcal H}^{+}_{g}\oplus {\mathcal 
H}^{-}_{g}$$
of the self-dual and anti-self-dual harmonic forms. This allows
us to uniquely write 
$$c_{1}=c_{1}^{+}+c_{1}^{-},$$
where the cohomology classes $c_{1}^{\pm}\in H^{2}(M,{\mathbb R})$ 
have harmonic representatives in ${\mathcal H}_{g}^{\pm}$, 
respectively.

For any given self-dual form $\eta$ on $(M,g)$, the 
corresponding {\em perturbed Seiberg-Witten equations}
\cite{taubes} read
 \begin{eqnarray} D_{A}\Phi &=&0 \label{dir} \\
 iF^+_A+\sigma (\Phi ) &=& \eta , \label{sd}\end{eqnarray}
where the unknowns are a Hermitian connection $A$ on 
the line bundle $L$ associated with $\mathfrak c$, and 
a section $\Phi$ of the twisted spinor bundle ${\mathbb V}_{+}$.
Here 
$D_{A}: \Gamma ({\mathbb V}_{+})\to \Gamma ({\mathbb V}_{-})$ 
is the Dirac operator determined by $A$, and 
$F_{A}^{+}$ is the  self-dual part of the curvature
of $A$, whereas $\sigma : {\mathbb V}_{+}\to \Lambda^{+}$
is the natural real-quadratic map induced by 
the isomorphism $\Lambda^{+}\otimes {\mathbb C} =
\odot^{2}{\mathbb S}_{+}$, with the conventional normalization  that
$|\sigma(\Phi )|^{2}= |\Phi |^{4}/8$.

Let ${\mathcal U}$ denote the affine space of differentiable 
unitary  connections
on $L$, and let  ${\mathcal V}$ denote the vector space of
differentiable sections of ${\mathbb V}_{+}$.  
The solution space 
$${\mathcal S}_{{\mathfrak c},g,\eta}=
\{ ~(A, \Phi )~|~ \mbox{(\ref{dir}) and (\ref{sd}) are satisfied} ~\}$$
of the perturbed Seiberg-Witten equations is thus a subset of 
${\mathcal U}\times {\mathcal V}$. Moreover, as long as the 
harmonic part of  $\eta$ is different from 
$2\pi c_{1}^{+}$, any solution of (\ref{dir}--\ref{sd})
will be  {\em irreducible}, in the sense that 
$\Phi \not\equiv 0$; thus, for generic $\eta$, 
$${\mathcal S}_{{\mathfrak c},g,\eta}\subset 
 {\mathcal U}\times ({\mathcal V}-0).$$ Now the 
{\em gauge  group} ${\mathcal G} = \{ u: M\to S^{1}\}$ acts on 
 ${\mathcal U}\times ({\mathcal V}-0)$ 
 by $(A,\Phi) \mapsto ( A - 2d\log u,  u\Phi)$, and 
 this action preserves the solution space ${\mathcal S}_{{\mathfrak 
 c},g,\eta}$. The quotient space 
 $${\mathcal M}_{{\mathfrak c},g,\eta}=
{\mathcal S}_{{\mathfrak c},g,\eta}/{\mathcal G}$$
 is called the Seiberg-Witten 
{\em moduli space} associated with the given 
spin$^{c}$ structure, and is tautologically
a subset of the  {\em configuation space} 
  $${\mathcal B}= [{\mathcal U}\times ({\mathcal V}-0)]/{\mathcal G},$$
  which is homotopy equivalent to 
  $T^{b_{1}(M)}\times {\mathbb C \mathbb P}_{\infty}$. 
For generic $\eta$, the moduli space is a smooth, compact
manifold of dimension 
$$\ell = \frac{c_{1}^{2}-(2\chi + 3\tau )}{4},$$
where $\chi$ and $\tau$ respectively denote the 
Euler characteristic and 
signature of $M$; in particular, 
if this integer is negative, the 
moduli space is empty for generic $\eta$. 
Moreover,  an orientation of the vector space 
$H^{1}(M,{\mathbb R})\oplus {\mathcal H}^{+}_{g}$
 determines an orientation of the moduli space.
Thus the homology class of 
${\mathcal M}_{{\mathfrak c},g,\eta}\subset {\mathcal B}$
gives us an element of 
$H_{\ell}({\mathcal B}, {\mathbb Z})\cong
H_{\ell}(T^{b_{1}(M)}\times {\mathbb C \mathbb P}_{\infty}, {\mathbb Z})$
which  turns out to be  independent of the metric $g$ and the
generic self-dual form $\eta$, and which  is  
 called \cite{ozsz,taubes3} the (generalized)
{\em	Seiberg-Witten invariant} of $(M,{\mathfrak c})$.  
 For our purposes, the only important point is that 
  when this   invariant is non-zero,
  the equations 
\begin{eqnarray} D_{A}\Phi &=&0 \label{specdir}\\
 -iF^+_A&=& \sigma (\Phi ) - t\phi , \label{specsd} \end{eqnarray}
must have a solution 
for any metric $g$, any self-dual harmonic $2$-form $\phi$, and any 
real constant $t$. Moreover, provided $\phi\not\equiv 0$, 
 this solution will be 
 irreducible 
for any sufficiently large $t$.

 \begin{lem} \label{erst}
Let $(M,g)$ be a  compact oriented Riemannian 4-manifold,
and let $\phi$ be a self-dual harmonic 2-form on
$(M,g)$. If, for a given spin$^{c}$ structure $\mathfrak c$,  the 
 perturbed Seiberg-Witten equations 
(\ref{specdir}--\ref{specsd}) have a solution for a given real
number $t$, then 
$$V^{1/3}\left( \int_{M}\left|\frac{2}{3}s+2w-t 
2\sqrt{2}| \phi |\right|^{3}d\mu 
 \right)^{2/3}\geq 8t^2 [\phi ]^2 - 32\pi t c_1 \cdot [\phi ],$$
 where $V$ denotes the total volume of $(M,g)$, 
 $w: M\to (-\infty, 0]$ is the lowest eigenvalue of 
 the self-dual Weyl curvature $W_{+}: \Lambda^{+}\to \Lambda^{+}$ of $g$, 
 $d\mu$ is the 
volume form of $g$, 
  $|\cdot |$ is the 
  point-wise norm     determined by  $g$,
and $c_1=c_1(L)$ is the first Chern class of 
the spin$^c$ structure $\mathfrak{c}$. 
\end{lem}

\begin{proof}
The Dirac equation $D_A\Phi =0$ implies the   Weitzenb\"ock formula
$$0= \frac{1}{2} \Delta |\Phi |^{2} + | \nabla \Phi|^{2}+
 \frac{s}{4} |\Phi |^2 + 2 \langle -iF^+_A , \sigma (\Phi )\rangle
.$$
For a solution of the perturbed Seiberg-Witten equations
(\ref{specdir}--\ref{specsd}), 
 we also have $-iF^+_A = \sigma (\Phi ) - t\phi$, so it follows that
 \begin{eqnarray} 
0 &=&  2 \Delta |\Phi |^{2} + 4| \nabla \Phi|^{2}+
 s |\Phi |^2 + |\Phi |^{4}-8t \langle \phi  , \sigma (\Phi )\rangle
 \nonumber \\
 &\geq& 2 \Delta |\Phi |^{2} + 4| \nabla \Phi|^{2}+
 s |\Phi |^2 + |\Phi |^{4}-t 2\sqrt{2}| \phi |~ | \Phi |^{2}.
 \label{wbk}
 \end{eqnarray}
 Multiplying by $|\Phi|^{2}$, we thus have
$$0\geq   2 |\Phi |^{2}\Delta |\Phi |^{2} + 4|\Phi |^{2}| \nabla \Phi|^{2}+
 (s -t 2\sqrt{2}| \phi |) |\Phi |^4 + |\Phi |^{6},
 $$
 and, upon integrating, we  obtain
$$
0 \geq  \int_{M}\left[ 2\left|d|\Phi |^{2}\right|^{2}+
4|\Phi |^{2}| \nabla \Phi|^{2}+
(s -t 2\sqrt{2}| \phi |) |\Phi |^4 + |\Phi |^{6}
\right]d\mu ,$$
so that 
\begin{equation}
\int \left[(-s) |\Phi |^4 - 4|\Phi |^{2}| \nabla \Phi|^{2} \right]d\mu
\geq \int \left[ |\Phi |^{6} -t 2\sqrt{2}| \phi |) |\Phi |^4 
\right]d\mu .
	\label{two}
\end{equation}

 Now recall that  any self-dual  2-form $\psi$ on any 
 oriented 4-manifold satisfies 
 the 
Weitzenb\"ock formula \cite{bourg}
$$(d+d^{*})^{2}\psi = \nabla^{*}\nabla \psi - 2W_{+}(\psi , 
\cdot ) + \frac{s}{3} \psi,$$
where $W_{+}$ is the self-dual Weyl tensor. 
It follows that 
$$
\int_{M}[-2W_{+}(\psi , \psi )] d\mu \geq 
\int_{M}(-\frac{s}{3})|\psi |^{2}~d\mu -
 \int_{M} |\nabla \psi |^{2}
 ~d\mu , 	
$$
so that 
 $$
  -\int_{M}2w|\psi |^{2} d\mu \geq 
\int_{M}(-\frac{s}{3})|\psi |^{2}~d\mu -
\int_{M} |\nabla \psi |^{2}
 ~d\mu , 
  $$ 
  and hence
   $$
  -\int_{M}(\frac{2}{3}s+2w)|\psi |^{2} d\mu \geq 
\int_{M}(-s)|\psi |^{2}~d\mu -
\int_{M} |\nabla \psi |^{2}
 ~d\mu . 
  $$ 
On the other hand, the  particular self-dual 2-form $\varphi = \sigma (\Phi 
)$ satisfies 
\begin{eqnarray*}
	|\varphi |^{2} & = & \frac{1}{8}|\Phi |^{4}  ,\\
	|\nabla \varphi |^{2} & \leq  & \frac{1}{2} |\Phi |^{2}|\nabla \Phi |^{2} .
\end{eqnarray*}
 Setting $\psi = \varphi$,  
  we thus have
  $$
   -\int_{M}(\frac{2}{3}s+2w)|\Phi |^{4} d\mu \geq 
\int_{M}\left[ (- s)|\Phi |^{4}-
 4 |\Phi |^{2} |\nabla \Phi |^{2}\right] d\mu .
  $$
 Combining this with (\ref{two}), we thus obtain  
 $$
 -\int_{M}(\frac{2}{3}s+2w)|\Phi |^{4} d\mu \geq  
 \int_{M} \left[ |\Phi |^{6} -t 2\sqrt{2}| \phi | |\Phi |^4 
\right]d\mu ,$$
and hence 
\begin{equation}
 -\int_{M}(\frac{2}{3}s+2w-t 2\sqrt{2}| \phi |)|\Phi |^{4} d\mu \geq  
 \int_{M}  |\Phi |^{6} .	
	\label{three}
\end{equation} 
By the H\"older inequality, this implies 
$$\left( \int \left|\frac{2}{3}s+2w
-t 2\sqrt{2}| \phi |\right|^{3}d\mu \right)^{1/3}
\left( \int |\Phi |^{6} d\mu \right)^{2/3}
\geq  
 \int|\Phi |^{6}~d\mu ,$$
 and hence  that 
 $$\int \left|\frac{2}{3}s+2w-t 2\sqrt{2}| \phi |\right|^{3}d\mu 
\geq  
 \int|\Phi |^{6}~d\mu .$$
 But  the H\"older inequality also tells us 
 that
 $$
V^{1/3} \left(\int|\Phi |^{6}~d\mu  \right)^{2/3}\geq \int |\Phi |^{4}d\mu ,
 $$
 so we now have 
\begin{eqnarray*}
V^{1/3}\left( \int_{M}\left|\frac{2}{3}s+2w-t 2\sqrt{2}| \phi |\right|^{3}d\mu 
 \right)^{2/3}&\geq&
\int    |\Phi |^4 d\mu\\  &=&8\int    |\sigma(\Phi) |^2 d\mu
\\&=& 8\int   |-iF^+_A+t\phi |^2 d\mu
\\&=& 8\int   \left( t^2|\phi |^2   -2t\langle iF^+_A,\phi
\rangle 
+  |iF^+_A|^2\right) d\mu
\\&= &  8\left(t^2  [\phi ]^2 - 2 t (2\pi c_1) \cdot [\phi ]
+  \int |iF^+_A|^2d\mu  \right)
\\&\geq& 8t^2 [\phi ]^2 - 32\pi t c_1 \cdot [\phi ]  ,
\end{eqnarray*}
as claimed.
\end{proof}

 \begin{lem} \label{prev}
 Let $\gamma =[g]$ be  a smooth conformal class
 on a smooth compact oriented 4-manifold $M$,
 and let $\phi \not\equiv 0$ be a closed $2$-form which is 
 self-dual with respect to $\gamma$. Suppose, 
 moreover, that for a fixed  spin$^{c}$ structure 
 $\mathfrak c$ and a fixed positive real number $t$ that
  the 
 perturbed Seiberg-Witten equations 
(\ref{specdir}--\ref{specsd}) have an irreducible solution for every metric
in the conformal class $\gamma$. Then, for any metric $g\in \gamma$,
 the scalar curvature $s$ and 
Weyl curvature $W_{+}$ satisfy the inequality 
$$ \int_{M}\left(\frac{2}{3}s+2w-t 
2\sqrt{2}| \phi |\right)^{2}d\mu 
 \geq 8t^2 [\phi ]^2 - 32\pi t c_1 \cdot [\phi ].$$
 \end{lem}
 \begin{proof}
 The key step in the  argument is a conformal rescaling trick,
 the general idea of which is due to Gursky \cite{G1}. 
 
 We begin by observing that there is a $C^{2}$ metric $g_{\gamma}\in \gamma$ 
 for which the function ${\mathfrak S}=s+3w-t 
3\sqrt{2}| \phi |$ is constant. Indeed, if $\hat{g}=u^{2}g$, 
the corresponding curvature quantity is given by 
$${\mathfrak S}_{\hat{g}} = u^{-3}\left( 6 \Delta_{g}u + {\mathfrak 
 S}_{g}u \right)$$
 because $w_{\hat{g}}=u^{-2}w_{g}$ and $|\phi|_{\hat{g}}
 =u^{-2}|\phi|_{g}$. The metric $g_{\gamma}$ 
 may therefore be constructed by minimizing the functional
 $${\mathcal F}(g) = 
 \frac{
 \int_{M} {\mathfrak S}_{g} d\mu_{g}
 }{
 \sqrt{\int_{M}~d\mu_{g}}
 }
$$
among metrics in the conformal class $\gamma$. The infimum of this functional
is negative because the Weitzenb\"ock formula (\ref{wbk}) shows  that every
metric in $\gamma$ has $s- t2\sqrt{2}| \phi |\leq 0$ somewhere, and
${\mathfrak S}\leq s- t2\sqrt{2}| \phi |$ everywhere, with strict inequality
 at any point where
$\phi \neq 0$.  Trudinger's approach to the Yamabe problem
thus produces a minimizer $g_{\gamma}$ of regularity $C^{2,\alpha}$ 
for any $\alpha < 1$.

Since  ${\mathfrak S}_{g_{\gamma}}$ is 
automatically a negative constant by the Euler-Lagrange equations,
we have 
$$\int_{M}\left(\frac{2}{3}s_{g_{\gamma}}+2w_{g_{\gamma}}-t 
2\sqrt{2}| \phi |_{g_{\gamma}}
 \right)^{2}d\mu_{g_{\gamma}} =
V^{1/3}_{g_{\gamma}}\left( \int_{M}\left|(\frac{2}{3}s_{g_{\gamma}}+ 
{2}w_{g_{\gamma}}-t 
2\sqrt{2}| \phi |_{g_{\gamma}} \right|^{3}d\mu_{g_{\gamma}} 
 \right)^{2/3},$$
 so that
 $$
 \int_{M}\left(\frac{2}{3}s_{g_{\gamma}}+2w_{g_{\gamma}}-t 
2\sqrt{2} | \phi |_{g_{\gamma}}
 \right)^{2}d\mu_{g_{\gamma}} \geq 8t^2 
 [\phi ]^2 - 32\pi t c_1 \cdot [\phi ]$$
 by the previous lemma. Thus the desired inequality at
 least holds for the particular metric $g_{\gamma}\in \gamma$.

We now compare the left-hand side with analogous 
expression for 
 the given metric $g$, following an idea of 
\cite{bcg1}. To do so, 
we express $g$ in the form $g=u^{2}g_{\gamma}$,
where $u$ is a positive $C^{2}$ function, and observe that 
\begin{eqnarray*}
\int_{M} {\mathfrak S}_{g} u^{2}d\mu_{g_{\gamma}}  
	 & = & \int u^{-3}\left(6\Delta_{g_{\gamma}}u+
	 {\mathfrak S}_{g_{\gamma}}u\right) u^{2}d\mu_{g_{\gamma}}
	 \\&=& \int \left(-6 u^{-2} |du|^{2}_{g_{\gamma}}+
	 {\mathfrak S}_{g_{\gamma}}\right)  d\mu_{g_{\gamma}}
	 \\ &\leq & \int 
	 {\mathfrak S}_{g_{\gamma}}  d\mu_{g_{\gamma}} .
\end{eqnarray*}
Applying Cauchy-Schwarz, we thus have
\begin{eqnarray*}
 -V^{1/2}_{g_{\gamma}}\left[\int {\mathfrak S}_{g}^{2} d\mu_{g}\right]^{1/2}&=&
	 -V^{1/2}_{g_{\gamma}}\left(\int 
	 {\mathfrak S}_{g}^{2} u^{4}d\mu_{g_{\gamma}}\right)^{1/2}\\
 & \leq  & 
 \int_{M} {\mathfrak S}_{g} u^{2}d\mu_{g_{\gamma}}
 \\&\leq &
 \int_{M} {\mathfrak S}_{g_{\gamma}}
	 d\mu_{g_{\gamma}}  \\
	 & = &  -
	 V^{1/2}_{g_{\gamma}}\left[\int_{M} 
	 {\mathfrak S}_{g_{\gamma}}^{2}d\mu_{g_{\gamma}}\right]^{1/2} . 
\end{eqnarray*}
Thus   
\begin{eqnarray*}\int_{M}\left(\frac{2}{3}s_{g}+2w_{g}-t 
2\sqrt{2} | \phi |_{g} \right)^{2}d\mu_{g} 
&=& \frac{4}{9}\int_{M}{\mathfrak S}_{g}^{2} d\mu_{g}\\
&\geq &
\frac{4}{9}\int_{M} 
	 {\mathfrak S}_{g_{\gamma}}^{2}d\mu_{g_{\gamma}}
	 \\&=&
\int_{M}\left(\frac{2}{3}s_{g_{\gamma}}+2w_{g_{\gamma}}-t 
2\sqrt{2} | \phi |_{g_{\gamma}}
 \right)^{2}d\mu_{g_{\gamma}} \\
&\geq& 
 8t^2 
 [\phi ]^2 - 32\pi t c_1 \cdot [\phi ]  ,
\end{eqnarray*}
 as claimed. 
 \end{proof}

 \begin{thm} \label{gest}
Let $M^4$ be a smooth compact oriented 4-manifold with
$b^{+}\geq 2$, and suppose that $\mathfrak{c}$
is a spin$^{c}$ structure  with  non-trivial 
  Seiberg-Witten invariant. 
 Let $g$ be any 
Riemannian metric on $M$, and let $\phi$ be a
$g$-self-dual harmonic 2-form with de~Rham class $[\phi]\in 
H^2(M,{\Bbb R})$. Then the scalar curvature $s$ and 
lowest eigenvalue $w$ of the self-dual Weyl curvature
$W_{+}$ of $g$
satisfy
$$\int \left(\frac{2}{3}s+2w\right) 
\frac{|\phi |}{\sqrt{2}} d\mu \leq 4\pi c_1\cdot [\phi
].$$
Here $d\mu$  and  $|\cdot |$ are respectively the 
volume form and 
  point-wise norm     determined by the metric $g$,
while $c_1=c_1(V_+)$ is the first Chern class of 
the spin$^c$ structure $\mathfrak{c}$. 
\end{thm}
\begin{proof}
 By Lemma \ref{prev}, we have 
$$
\int_{M}\left(\frac{2}{3}s+2w-t 2\sqrt{2}| \phi |\right)^{2}d\mu	
\geq 8t^2 [\phi ]^2 - 32\pi t c_1 \cdot [\phi ]
$$
for all $t>0$. Thus
$$8t^2 [\phi ]^2-4\sqrt{2}t \int_{M}(\frac{2}{3}s+2w) |\phi| d\mu
+\int_{M}(\frac{2}{3}s+2w)^{2}d\mu \geq 
8t^2 [\phi ]^2 - 32\pi t c_1 \cdot [\phi ],$$
and hence 
$$\int_{M}(\frac{2}{3}s+2w) \frac{|\phi|}{\sqrt{2}} d\mu - \frac{8}{t}
\int_{M}(\frac{2}{3}s+2w)^{2}d\mu\leq 4\pi  c_1 \cdot [\phi ]. $$
Taking the limit as $t\to \infty$ 
 then yields the desired result.
\end{proof}

\begin{cor}
Let $M^4$ be a smooth compact oriented 4-manifold with
$b^{+}\geq 2$, and suppose that $\mathfrak{c}$
is a spin$^{c}$ structure  with  non-trivial  
Seiberg-Witten invariant. 
 Let $g$ be any 
Riemannian metric on $M$, and let $\phi$ be a
$g$-self-dual harmonic 2-form with de~Rham class $[\phi]\in 
H^2(M,{\Bbb R})$. Then the scalar curvature $s$ and Weyl curvature
$W_{+}$ of $g$
satisfy
$$\int \left(s-\sqrt{6}|W_{+}|\right) 
|\phi | d\mu \leq 6\sqrt{2}\pi c_1\cdot [\phi
].$$
\end{cor}
\begin{proof} Because $W_{+}$ is a  trace-free endomorphism of 
$\Lambda^{+}$, 
$$-\sqrt{\frac{2}{3}}|W_{+}|\leq w.$$
Substituting this into Theorem \ref{gest},  and multiplying by 
$3/\sqrt{2}$, 
we thus obtain the desired result.
\end{proof}

Since Theorem \ref{gest} applies 
  to  every metric conformal
to a given $g$, we can improve it, as follows:

 \begin{thm} \label{trigger}
Let $M^4$ be a smooth compact oriented 4-manifold with
$b^{+}\geq 2$, and suppose that $\mathfrak{c}$
is a spin$^{c}$ structure  with  non-trivial  
Seiberg-Witten invariant. 
 Let $g$ be any 
Riemannian metric on $M$, and let $\phi$ be a
$g$-self-dual harmonic 2-form with de~Rham class $[\phi]\in 
H^2(M,{\Bbb R})$. Then 
the function  $f = \sqrt{|\phi |}$
 satisfies  
$$\int_{M} \left(\frac{2}{3}s_{g}+2w_{g}\right) 
|\phi |_{g} d\mu_{g} 
+4\int_{M} |df|_{g}^{2}d\mu_{g}  
\leq (4\pi \sqrt{2}) c_1\cdot [\phi
].$$
\end{thm}
\begin{proof} 
We may obviously assume that $\phi\not\equiv 0$, since
otherwise there is nothing to prove. 

Now observe that, 
for any smooth positive function $u$ on $M$,
 the metric $\hat{g}=u^{2}g$
satisfies 
$$(\frac{2}{3}s_{\hat{g}}+2w_{\hat{g}}) 
|\phi |_{\hat{g}} ~ d\mu_{\hat{g}}=
  u^{-1}|\phi|_{g}\left[4\Delta_{g}u+(\frac{2}{3}s_{g}+2w_{g})u\right]
d\mu_{g},$$
and, since $|\phi|$ has bounded derivative, we may  integrate by
parts to obtain 
$$\int_{M} \left(\frac{2}{3}s_{\hat{g}}+2w_{\hat{g}}\right) 
|\phi |_{\hat{g}} ~d\mu_{\hat{g}} =
\int_{M} \left(\frac{2}{3}s_{g} +2w_{g} \right)
|\phi |_{g} ~d\mu_{g} 
+4\int_{M} \langle d (u^{-1} |\phi |), du \rangle_{g} d\mu_{g}.$$
Applying Theorem \ref{gest} to 
$\hat{g}$ thus gives us 
\begin{equation}
\int_{M} \left(\frac{2}{3}s_{g}+2w_{g}\right) 
|\phi |_{g} d\mu_{g} 
+4\int_{M} \langle d (u^{-1} |\phi |), du \rangle_{g} d\mu_{g}
\leq 	(4\pi \sqrt{2}) c_1\cdot [\phi].
	\label{roy}
\end{equation}

Now, for some  
$\varepsilon > 0$,  let us take $u=u_{\varepsilon}$ to be given by
$$u_{\varepsilon}= \sqrt{\alpha_{\varepsilon}(|\phi|)}, $$
where $\alpha_{\varepsilon}: {\mathbb R}\to {\mathbb R}$ is
a smooth function with 
\begin{eqnarray*}
	\alpha_{\varepsilon}(x)  & = & x ~~\forall x \in 
		[\varepsilon,\infty),    
		\\
		\alpha_{\varepsilon}^{\prime}(x) & = & 0 ~~\forall  x \in 
	(-\infty, \frac{\varepsilon}{2}], \mbox{ and}\\
	\alpha_{\varepsilon}^{\prime\prime}(x) & \geq & 0 ~~\forall  x\in {\mathbb R}. 	
\end{eqnarray*}
Then, since  $\alpha_{\varepsilon}^{\prime}(x) \leq 1$ for all $x$, 
\begin{eqnarray*}
	\langle d(u_{\varepsilon}^{-1}|\phi |) , du_{\varepsilon} 
	\rangle_{g} & = & \frac{\alpha_{\varepsilon}^{\prime}(|\phi |)}{2\alpha_{\varepsilon}(|\phi |)}
\left(1-\frac{|\phi|\alpha_{\varepsilon}^{\prime}(|\phi |)}{2\alpha_{\varepsilon}(|\phi |)}\right)|~d|\phi|~|^{2}_{g}
  \\
	 & \geq & \frac{1}{4} 
\frac{\alpha_{\varepsilon}^{\prime}(|\phi |)}{\alpha_{\varepsilon}(|\phi |)}
~|~d|\phi|~|^{2}_{g},
\end{eqnarray*}
with equality when $|\phi | \geq \varepsilon$. 
Hence $\langle d(u_{\varepsilon}^{-1}|\phi |) , du_{\varepsilon} \rangle
\geq 0$ everywhere, and 
$\langle d(u_{\varepsilon}^{-1}|\phi | ), du_{\varepsilon} \rangle = 
|df|^{2}$ on the set  $M-M_{\varepsilon}$ where $|\phi | \geq 
\varepsilon$. Thus 
$$\int_{M} \langle d (u^{-1} |\phi |), du \rangle ~ d\mu ~ 
> ~ \int_{M-M_{\varepsilon}}|df|^{2}d\mu , $$
and (\ref{roy}) therefore implies
$$
\int_{M} \left(\frac{2}{3}s+2w\right) 
|\phi | ~d\mu 
+4\int_{M-M_{\varepsilon}} |df|^{2} d\mu
<	(4\pi \sqrt{2}) c_1\cdot [\phi].
$$
However, $\varepsilon > 0$ is arbitrary, and 
$\cap_{\varepsilon}M_{\varepsilon}$ is the nodal set of $\phi$, which
has measure zero. Taking the supremum of the left-hand side over 
$\varepsilon  > 0$
thus gives 
us the desired inequality.
\end{proof}

\section{Vanishing Theorems}

One of Witten's most elegant observations is the 
fact that a given $4$-manifold can only have finitely many
spin$^{c}$ structures for which the Seiberg-Witten invariant
is non-zero; this follows from the fact that there is an 
{\em a priori} upper bound on 
$(c_{1}^{+})^{2}$ in terms of  scalar curvature \cite{witten}. 
When the $4$-manifold admits a hyperbolic metric, however, 
we will now see that the invariant must also vanish for 
most of the  spin$^{c}$ structures which slip under the 
bar of  
Witten's upper bound.

\begin{main}\label{care}
Let $(M,g)$ be a compact, anti-self-dual Einstein manifold 
with scalar curvature
$s=-12$. 
Suppose, moreover, that $b^{+}(M)\geq 2$, and 
 that $c_{1}$ is the first Chern class of a
spin$^{c}$ structure on $M$ for which the (generalized) Seiberg-Witten
invariant is non-zero; let $c_{1}^{+}$ denote its
self-dual part of $c_{1}$ with respect to $g$. Then
$$\frac{(c_{1}^{+})^{2}}{2\chi+3\tau}\leq \frac{(2-\nu 
)^{2}\cos^{2}\theta}{3} ,$$
where $\chi$ and $\tau$ respectively denote the 
Euler characteristic and signature of $M$, and where 
 the invariants   $\theta$ and $\nu$ of $(M,g)$ are as in  
 Definitions \ref{angl} and \ref{quot}.  
\end{main}
\begin{proof}
For any self-dual harmonic $2$-form on $(M,g)$, Theorem \ref{trigger}
tell us that 
$$\int_{M} \frac{2}{3}s_{g}
|\phi |_{g} d\mu_{g} 
+4\int_{M} |d\sqrt{|\phi |}|_{g}^{2}d\mu_{g}  
\leq (4\pi \sqrt{2}) c_1\cdot [\phi
]$$
because the anti-self-duality of $g$ guarantees that 
 $w_{g}\equiv 0$. Since $s$ is constant, the definition
 of $\nu (M,g)$ thus tells us that 
 $$\frac{2}{3}s +4\nu \leq \frac{\int_{M} \frac{2}{3}s_{g}
|\phi |_{g} d\mu_{g} 
+4\int_{M} |d\sqrt{|\phi |}|_{g}^{2}d\mu_{g} }{\int |\phi |d\mu_{g}} \leq 
 (4\pi \sqrt{2}) \frac{c_1\cdot [\phi]}{\int |\phi |d\mu_{g}}.$$
 Taking $[\phi ] =-c_{1}^{+}$, we therefore have
 \begin{eqnarray*}
  	-s-6\nu & \geq  &   6\pi \sqrt{2} 
 \frac{(c_{1}^{+})^{2}}{\int |\phi |d\mu_{g}} \\
  	 & = &  6\pi \sqrt{2} 
 \frac{(\int |\phi |^{2}d\mu_{g})^{1/2}}{\int |\phi |d\mu_{g}} 
 \sqrt{(c_{1}^{+})^{2}}  \\
  	 & \geq & 6\pi \sqrt{2}\frac{V^{-1/2}}{\cos \theta } q 
  	 \sqrt{2\chi + 3\tau}
  \end{eqnarray*} 
 where
 $$q=\sqrt{\frac{(c_{1}^{+})^{2}}{2\chi + 3\tau}}.$$
 However, we have the Gauss-Bonnet-like formula \cite{hit} 
 $$2\chi + 3\tau  
=\frac{1}{4\pi^{2}}\int_{M}\left(
2|W_{+}|^{2}+\frac{s^{2}}{24} -\frac{|\stackrel{\circ}{r}|^{2}}{2}
\right) d\mu ,	$$
 for any metric on $M$, and for our anti-self-dual Einstein metric
 $g$ this simplifies to become
 $$2\chi + 3\tau = \frac{s^{2}V}{96\pi^{2}}.
 $$
 Thus, since  $s< 0$, 
 $$
 |s|-6\nu \geq \sqrt{72\pi^{2}}
 \frac{ V^{-1/2}}{\cos \theta}q
 \sqrt{\frac{s^{2}V}{96\pi^{2}}}= \frac{\sqrt{3}}{2}\frac{|s| q}{\cos 
 \theta},
 $$
and 
 $$ \frac{1}{\sqrt{3}}(2-\frac{12}{|s|}\nu )\cos \theta\geq q .$$
 With the normalization $s=-12$, this then gives us 
 $$  
 \frac{1}{3} (2-\nu )^{2}\cos^{2}\theta 
 \geq q^{2}=\frac{(c_{1}^{+})^{2}}{2\chi + 3\tau},$$
 as claimed. 
\end{proof}

\begin{cor}
Let $(M,g)$ be a compact, anti-self-dual Einstein manifold 
with scalar curvature
$s=-12$ and $b^{+}(M)\geq 2$. Then 
the Seiberg-Witten invariant vanishes for any 
spin$^{c}$ structure for which 
$$(c_{1}^{+})^{2} > \frac{(2-\nu 
)^{2}\cos^{2}\theta}{3} (2\chi + 3\tau)  
$$
or for which 
 $$
|(c_{1}^{-})^{2}| > \frac{(2-\nu 
)^{2}\cos^{2}\theta-3}{3} (2\chi + 3\tau) 
$$
\end{cor}
\begin{proof} 
If the Seiberg-Witten invariant were non-zero, we would necessarily 
have
$$(c_{1}^{+})^{2}-|(c_{1}^{-})^{2}| = c_{1}^{2} \geq 2\chi + 3\tau$$
because the virtual dimension of the moduli space must be 
non-negative. Theorem \ref{care} thus guarantees that 
$$\left( \frac{(2-\nu 
)^{2}\cos^{2}\theta}{3} -1\right) (2\chi + 3\tau) \geq 
(c_{1}^{+})^{2}-(2\chi + 3\tau)\geq |(c_{1}^{-})^{2}|.$$
This proves the corollary  by contraposition. 
\end{proof}

In particular, since $|(c_{1}^{-})^{2}|\geq 0$ for all spin$^{c}$
structures,  we obtain 

\begin{cor}
Let $(M,g)$ be a compact, anti-self-dual Einstein manifold 
with scalar curvature
$s=-12$ and $b^{+}(M)\geq 2$. If 
$$\nu > 
2- 
\sqrt{3} \sec \theta  , 
$$
then 
all the Seiberg-Witten invariants of $M$ must vanish for
the given orientation. 
\end{cor}

Theorem \ref{wild} now follows immediately, since 
 $\theta\neq 0$  for the  spaces in question. 
 But it remains to be seen, of course, whether this result 
 is actually non-vacuous! Can one at least show
 that $\nu \geq 2-\sqrt{3}$
  for {\em some}
 real-hyperbolic $4$-manifolds?

\bigskip 
\noindent
{\bf Acknowledgment} ~~The author would like to thank Cliff Taubes for
some helpful comments, and the Max-Planck-Institut f\"ur 
Mathematik in den Naturwissenschaften, Leipzig, for
its hospitality during the inception  of this work.

  \end{document}